\documentclass[11pt]{article}
\usepackage[english,ukrainian]{babel}
\usepackage{amssymb,amsfonts}
\usepackage{enumerate}
\usepackage{amsmath,amsthm}
\def\bey#1\eey{\begin{equation}#1\end{equation}}
\def\ben#1\een{\begin{equation*}#1\end{equation*}}

\newcommand{\rom}[1]{\textrm{\mdseries\upshape#1}}

\newcommand{\supp}{\mbox{supp\,}}
\newcommand{\dist}{\mbox{dist\,}}
\newcommand{\ov}{\overline}
\newcommand{\Z}{{\Bbb Z}}
\newcommand{\R}{{\Bbb R}}
\newcommand{\N}{{\Bbb N}}
\newcommand{\D}{{\Bbb D}}
\renewcommand{\C}{{\Bbb C}}
\newcommand{\Cvar}{{\Bbb C}}
\newcommand{\ve}{\varepsilon}
\newcommand{\ity}{\infty}

\renewcommand{\th}{\theta}

\newcommand{\de}{\delta}
\newcommand{\De}{\Delta}

\newcommand{\al}{\alpha}
\newcommand{\vfi}{\varphi}
\newcommand{\si}{\sigma}

\newcommand{\Bl}{\Bigl(}
\newcommand{\Br}{\Bigr)}
\newcommand{\Bm}{\Bigr|}

\renewcommand{\Im}{\mathop{\rm Im}}
\renewcommand{\Re}{\mathop{\rm Re}}
\newcommand{\SH}{\mathop{\rm SH}}

\newcommand{\ze}[1]{\zeta_l^{(#1)}}
\newcommand{\xxi}[1]{\xi_l^{(#1)}}
\newcommand{\pl}{P^{(l)}}
\newcommand{\ql}{Q^{(l)}}

\newcommand{\ml}{\mu^{(l)}(\zeta)}
\renewcommand{\L}{{\cal L}}

\newcommand{\diam}{\mathop{\rm diam}}

\newtheorem{theorem}{Theorem}
\newtheorem{corollary}{Corollary}
\newtheorem{lemma}{Lemma}

\newtheorem{pheo}{Theorem}

\numberwithin{equation}{section}
\begin{document}


{\bf Approximation of subharmonic functions in the unit disk}

\medskip
{\bf I.E. Chyzhykov}


\bigskip

Chyzhykov Igor Elbertovych

Faculty of Mechanics and Mathematics,

Lviv Ivan Franko National University

Universytets'ka 1, Lviv, 79000

ichyzh@lviv.farlep.net

\date{}

\section{Introduction}
We use the standard notions of  subharmonic function theory
\cite{Hay}. Let us introduce some notation. Let $U(E,t)=\{\zeta\in
\Cvar: \dist(\zeta, E)<t\}$, $E\subset \C$,   $t>0$, where $\dist
(z, E)\stackrel {\rm def} = \inf_{\zeta \in E} |z-\zeta|$, and
$U(z,t)\equiv U(\{z\},t)$ for  $z\in \Cvar$.
 The class of
subharmonic functions in a domain $G\subset \C$ is denoted by
$\SH(G)$. For a subharmonic function $u\in \SH(U(0,R))$, $0<R\le
+\infty$ we write $B(r,u)=\max \{ u(z): |z|=r\}$, $0<r<R$, and
define the order $\rho[u]$ by \newline 
$\rho[u]=\limsup\limits_{r\to+\ity} \log B(r,u)/\log r$ if
$R=\infty$ and by $\si[u]=\limsup\limits_{r\to R} \log B(r,u)/\log
\frac1{R-r}$ if $R<\infty$.

 Let also $\mu_u$ denote the Riesz
measure associated with the subharmonic function $u$, $n(r, u)=
\mu_u(\overline{U(0,r)})$, $m$ be the planar Lebesgue measure, $l$
be the Lebesgue measure on the positive ray. For an analytic
function $f$ in $\mathbb{D}$ we write $Z_f=\{z\in \D: f(z)=0\}$.
The symbol $C(\cdot)$ with indices stands for some positive
constants depending only on values in the brackets. We write
$a\asymp b$ if $C_1 a\le b\le C_2 a$ for some positive constants
$C_1$ and $C_2$, and $a(r)\sim b(r)$ if $\lim _{r\to R}
a(r)/b(r)=1$.

An important result was proved by  R.\ S.~Yulmukhametov
\cite{Yul}. {\it For any function $u\in \SH(\Cvar)$ of order
$\rho\in (0,+\ity)$, and $\alpha>\rho$ there exist an entire
function $f$ and a set $E_\al\subset \Cvar$ such that
\begin{equation}
\big|u(z)-\log   |f(z)|\big|\le C(\al) \log  |z|, \quad  z\to\ity,
\; z\not\in E_\al, \label{e1.1}
\end{equation}
and $E_\al$ can be covered by a family of disks $U(z_j,t_j)$,
$j\in \N$, with $\sum_{|z_j|>R} t_j=O(R^{\rho-\al})$,
$(R\to+\ity)$.}

If $u\in \SH(\D)$ a counterpart of (\ref{e1.1}) holds with $\log
\frac 1{1-|z|}$ instead  of $\log|z|$ and an apropriate choice of
$E_\al$.


 From the recent result of Yu.\
Lyubarskii and Eu.\ Malinnikova \cite{LM} it follows that for
$L_1$  approximation relative to planar measure, we may drop the
assumption that $u$ has  finite order of growth, and obtain
sharp estimates.

\begin{pheo}[\cite{LM}] 
\label{th:in2} Let  $u\in \SH(\Cvar)$.
Then, for each $q>1/2$, there exist
  $R_0>0$ and an entire function  $f$ such
that
\begin{equation}
\frac 1{\pi R^2} \int\limits_{|z|<R}
         \big|u(z)-\log |f(z)|\big|dm(z) <q \log  R, \quad R> R_0.
\label{eq:in1e}
\end{equation}
\end{pheo}
An example constructed in \cite{LM} shows that we cannot  take
$q<1/2$ in estimate (\ref{eq:in1e}). The case $q=1/2$ remains open.


The following theorem complements this result.

Let $\Phi$  be the class of slowly growing functions $\psi\colon
[1,+\ity)\to (1,+\ity)$ (in particular, $\psi(2r)\sim \psi(r)$ as
$r\to+\ity$).
\begin{pheo}[\cite{AA}]\label{te:1}
Let $u\in \SH({\Bbb C})$, $\mu=\mu_u$. If for
some $\psi\in \Phi$ there exists a constant $R_1$ satisfying the
condition
\begin{equation}
(\forall R>R_1): \mu(\{ z:R<|z|\le R\psi(R)\})>1,
\label{e1.3}
\end{equation}
then there exists an entire function $f$ such that $(R\ge R_1)$
\begin{gather}
\int\limits_{|z|<R}
         \big|u(z)-\log  |f(z)|\big|\,dm(z) = O(R^2 \log   \psi (R)).
\label{e1.4}
\end{gather}
\end{pheo}

\noindent{\it Remark 1.1.} In the case $\psi(r)\equiv q>1$ we obtain
Theorem~1~\cite{LM}.

The following example and Theorem \ref{te:2} show (see \cite{AA} for details) that estimate
(\ref{e1.4}) is sharp in the class of subharmonic functions
satisfying (\ref{e1.3}).

For $\vfi\in \Phi$, let
\begin{equation*} 
u(z)= u_\vfi(z)=\frac12 \sum_{k=1}^{+\ity} \log   \Bigl| 1-\frac
z{r_k}\Bigr|,
\end{equation*}
where $r_0=2$, $r_{k+1}=r_k \vfi(r_k)$, $k\in \N\cup \{ 0\}$.
Thus, $\mu_{u}$ satisfies condition (\ref{e1.3}) with
$\psi(x)=\vfi^3(x)$.
\begin{pheo} \label{te:2}
Let  $\psi\in\Phi$ be such that $\psi(r)\to +\infty$
$(r\to+\infty)$. There exists no entire function $f$ for which
\begin{equation*}
\int\limits_{|z|<R}
         \big|u_\psi(z)-\log  |f(z)|\big|\,dm(z) = o(R^2 \log   \psi (R)), \quad
R\to\infty.
\end{equation*}
\end{pheo}

A further question appears naturally: Are there counterparts of
Theorems \ref{th:in2} and \ref{te:1} for subharmonic functions in
the unit disk? We have the following theorem.

\begin{theorem}\label{t:1}
\it Let $u\in \SH(\D)$. There exist an absolute constant $C$ and
an analytic function $f$ in $\D$ such that \bey\int_\D
\bigl|u(z)-\log|f(z)|\bigr|\, dm(z)<C.  \label{e:11} \eey
\end{theorem}

For a measurable set $E\subset [0,1)$ we define the density
$$\mathcal{D}_1 E=\varlimsup_{R\uparrow 1} \frac{l(E\cap [R,1))}{1-R}.$$

\begin{corollary} Let $u\in \SH(\D)$, $\varepsilon >0$. There
exist an analytic function $f$ in $\D$ and $E\subset [0,1)$, $\mathcal{D}_1
E<\ve$, such that \bey \label{e:12} \int_0^{2\pi}
\bigl|u(re^{i\th})-\log|f(re^{i\th})|\bigr| d\th= O\Bigl( \frac
1{1-r}\Bigr), \quad r\uparrow 1, r\not \in E. \eey
\end{corollary}

Relationship (\ref{e:12}) is equivalent to the condition $$
T(r,u)-T(r,\log|f|)=O((1-r)^{-1}), \quad r\uparrow 1, r\not \in
E,$$ where $T(r,v)$ is the Nevanlinna characteristic of
a subharmonic function~$v$. The author does not know whether (\ref{e:12}) is best possible.

\noindent{\it Remark 1.2.} No restriction on the Riesz measure
$\mu_u$ or  the growth of $u$ is required in Theorem~\ref{t:1}.

\noindent {\it Remark 1.3.} It is clear that (\ref{e:11}) is sharp
in the  class
 $\SH(\D)$, but under  growth restrictions  can be
 improved.
 \begin{pheo}[Hirnyk {\cite{Hi}}] \sl Let $u\in \SH(\D)$, $\si[u]<+\infty$. Then  there
 exists an analytic function $f$ in $\D$ such that
\ben\int_0^{2\pi}  \bigl|u(re^{i\th})-\log|f(re^{i\th})|\bigr|
d\th= O\Bigl (\log ^2 \frac 1{1-r}\Bigr), \quad r\uparrow 1. \een
\end{pheo}
%

Theorem \ref{t:1} does not allow the conclusion that \bey
u(z)-\log|f(z)|=O(1), \quad z\in \D\setminus E \label{e:unio}\eey
for any ``small''\ set $E$.

Sufficient conditions for (\ref{e:unio}) in the complex plane were
obtained in \cite{LM}. It uses so called notion of a {\it locally
regular measure} which admits a {\it partition of slow variation.}

We also prove a counterpart of Theorem $3'$ of \cite{LM} using a
similar concept. A~corresponding Theorem~\ref{t:3} will be
formulated in section \ref{s:3}. Here we formulate an application
of Theorem~\ref{t:3}.

\begin{theorem} \label{t:2} \sl Let $\gamma_j=(z=z_j(t): t\in [0,1])$, $1\le
j \le m$ be smooth Jordan curves in $\overline{U(0,1)}$
 such that $\arg
z_j(t)=\theta_j(|z_j(t)|)\equiv \theta_j(r)$, $|z_j(1)|=1$,
$|\theta_j'(r)|\le K$ for $r_0\le r<1$ and some constants $r_0\in
(0,1)$, $K>0$, $1\le j\le m$.  Let $u\in SH(\D)$, $\supp
\mu_u\subset \bigcup_{j=1}^m [\gamma_j]$, $\mu_u([\gamma_j]\cap
[\gamma_k])=0$, $j\ne k$,  and
$$\mu_u\Bigl|_{[\gamma_j]}(U(0,r))=\frac{\Delta_j}{(1-r)^{\sigma(r)}}
,$$ where $\Delta_j$ is a positive constant,
$\sigma(r)=\rho\bigl(\frac 1{1-r}\bigr)$, $\rho(R)$ is a proximate
order \cite{Le}, $\rho(R)\to\sigma>0$ as $R\to+\infty$.

Then there exists an analytic function $f$ such that for all
$\varepsilon>0$
\begin{gather}
\log|f(z)|-u(z)=O(1), \label{e:uni1}
\end{gather} $ z\not\in E_\varepsilon=\{\zeta\in \D: \dist
(\zeta, Z_f)\le \varepsilon(1-|\zeta|)^{1+\sigma(r)}\}$, 
where
\begin{gather} \label{e:uni11}
 \log|f(z)|-u(z)\le C,
\end{gather}
for some $C>0$ and all $z\in \D$. Moreover, $$Z_f \subset \bigcup
_{\zeta \in \bigcup_j[\gamma_j]} U(\zeta,
2(1-|\zeta|)^{1+\sigma(r)}),$$ and
\begin{equation} \label{e:tob}
T(r,u)-T(r,f)=O(1), \quad r\uparrow 1.
\end{equation}
\end{theorem}

\noindent{\it Remark 1.4.} Obviously, we can't obtain a lower
estimate for the left-hand side of (\ref{e:uni11}) for all $z$,
because it   equals $-\infty$ on $Z_f$.

 Theorems similar to Theorem~\ref{t:2} are proved in \cite[Ch.10,Th.10.16,10.20]{Hay2}.
The difference is that in \cite{Hay2} only more crude estimates are obtained  for approximation in a more general settings.

\section{Proof of Theorem \ref{t:1}}
\subsection*{2.1. Preliminaries} Let $u\in \SH(\D)$. Then
the Riesz measure $\mu_u$ is finite on compact subsets of $\D$. In
order to apply a partition theorem (Theorem E) we have to modify
the Riesz measure. On subtracting  an integer-valued
discrete measure $\tilde \mu$ from $\mu_u$ we may arrange that
$\nu(\{p\})=(\mu_u-\tilde \mu)(\{p\})<1$ for any point $p\in \D$.
The measure $\tilde \mu$ corresponds to the zeros of an entire
function $g$. Thus we can consider $\tilde u=u-\log |g|$,
$\mu_{\tilde u}=\nu$. According to Lemma~1~\cite{AA} in any
neighbourhood of the origin there exists a point $z_0$ with the
following properties:
\begin{itemize}
\item[a)]  on each line $L_\alpha$ going through   $z_0
$
there is at most one point  $\zeta_\alpha$ such that
 $\nu(\{\zeta_\alpha\})>0$, while  $\nu(L_\al\setminus \{\zeta_\al\})=0;$
\item[b)]  on each circle $K_\rho$ with center
 $z_0$ there exists at most one point $\zeta_\rho$ such
that
 $\nu(\{\zeta_\rho\})>0$, while $\nu(K_\rho\setminus
\{\zeta_\rho\})=0$.
\end{itemize}

 As it follows from the proof of Lemma 1~\cite{AA}, the set of points $z_0$ that
do not satisfy one of the conditions a) and b) has planar measure
zero. The similar assertion holds for the polar set
$u(z_0)=-\infty$ \cite[Chap.5.9,Theorem 5.32]{Hay}. Therefore, we can assume  that
properties  a), b) hold, and $u(z_0)\ne -\infty$.

Then consider the subharmonic function $u_0(z)=u\Bigl( \frac
{z_0-z}{1-z\bar z_0}\Bigr)\equiv u(w(z))$, $u_0(0)=u(z_0)$. Since
$|w'(z)|= \frac {1-|z_0|^2}{|1-z\bar z_0|^2}$,  we have
$||w'(z)|-1|\le 3|z_0|$ for $|z_0| \le 1/2$.

The Jacobian of the transformation $w(z)$ is $|w'(z)|^2$,
consequently this change of variables doesn't change relation
(\ref{e:11}).

Let
\begin{equation}\label{e:u3}
 u_3(z)=\int_{U(0,1/2)} \log|z-\zeta|\, d\mu_u(\zeta).
\end{equation}

The subharmonic function $u(z)- u_3(z)$ is harmonic in $U(0,1/2)$.


 Let $q\in(0,1)$ be such that
\begin{equation} \sum_{j=1}^{12} q^{j} >11. \label{e:qcon}
\end{equation}
We define $(n\in\{0, 1, \dots\})$
 $$R_n=1-q^n/2,\; A_n=\{\zeta:
R_n\le |\zeta|<R_{n+1}\},\;  M_n=M_n(q)=\Bigl[ \frac{2\pi}{\log
\frac{R_{n+1}}{R_n}}\Bigr], \; $$ $$ A_{n,m}=\Bigl\{\zeta\in A_n:
\frac{2\pi m}{M_n} \le \arg_0\zeta <\frac{2\pi (m+1)}{M_n}\Bigr\},
\quad 0\le m\le M_n-1.$$

Represent $\mu_u\Bigl|_{A_{n,m}}=\mu_{n,m}^{(1)}+\mu_{n,m}^{(2)}$
such that
\begin{itemize}
  \item [i)] $\supp \mu_{n,m}^{(j)} \subset \ov{A}_{n,m}$, $j\in\{1,2\}$;
    \item [ii)] $\mu_{n,m}^{(1)}(\ov{A}_{n,m})\in 2\Z_+$, $0\le  \mu_{n,m}^{(2)}
(\ov{A}_{n,m})<2$.
\end{itemize}

Let $$\mu_{n}^{(j)}=\sum_{m=0}^{M_n-1} \mu_{n,m}^{(j)}, \quad
\tilde \mu^{(j)}=\sum_n \mu_n^{(j)}, \; j\in\{1,2\}.$$

Property ii) implies
\begin{equation}\label{e:esm2}
\mu_n^{(2)}(\ov{A_n})\le \frac{13}{(1-q)(1-R_n)},\quad
n\to+\infty,
\end{equation}
as follows from the asymptotic equality
\begin{equation} \log \frac{R_{n+1}}{{R_n}}\sim (1-q)(1-R_n),
\quad n\to+\infty, \label{e:22'}
\end{equation}
and the definition of $M_n$.

Let \bey u_2(z)=\int_{\D} \log\Bigl|E\Bigl(
\frac{1-|\zeta|^2}{1-\bar \zeta z}, 1\Bigr) \Bigr| d\tilde
\mu^{(2)}(\zeta), \label{e:can} \eey where
$E(w,p)=(1-w)\exp\{w+w^2/2+\dots+w^p/p\}$, $p\in \N$ is the
Weierstrass primary factor.

\begin{lemma} $u_2\in SH(\D)$, and
$$T(r,u_2)=O\Bl\log ^2 \frac 1{1-r}\Br, \; r \uparrow 1, \quad
\int_{\D} |u_2(z)|\, dm(z)<C_1(q).$$ \label{l:d1}
\end{lemma}
\begin{proof}[Proof of Lemma \ref{l:d1}]
The following estimates for $\log|E(w,p)|$ are well-known
(cf. \cite[Ch.1, \S 4, Lemma 2]{Le}, \cite[Ch.4.1, Lemma
4.2]{Hay})
\bey\label{e:pfe}
\begin{split}
|\log E(w,1)| \le \frac {|w|^2}{2(1-|w|)}, \quad |w|<1,\\
\log |E(w,1)|\le 6e |w|^2, \quad w\in \C.
\end{split}
\eey First, we prove convergence of the integral in (\ref{e:can}).
For fixed  $R_n$ let $|z|\le R_n$. We choose $p$ such that $q^p<
1/4$. Then for $|\zeta|\ge R_{n+p}$ we have $$
\frac{1-|\zeta|^2}{|1-\bar \zeta z|}\le \frac{2(1-|\zeta|)}{1-
|z|}\le \frac {2(1-R_{n+p})}{1-R_n} < \frac 12.$$ Hence, using the
first estimate (\ref{e:pfe}), (\ref{e:esm2}) and the definition of
$R_n$ we obtain
\begin{gather*}
\int_{|\zeta|\ge R_{n+p}} \Bigl|\log \Bigl| E\Bigl( \frac{1-|\zeta|^2}{1-\bar \zeta
z},1\Bigr)
\Bigr|\Bigr| \, d\tilde\mu^{(2)}(\zeta)\le
\int_{|\zeta|\ge R_{n+p}} \Bigl( \frac{2(1-|\zeta|)}{1-|z|}
\Bigr)^2
 \, d\tilde\mu^{(2)}(\zeta)\le\\ \le
 \frac {4}{(1-|z|)^2} \sum_{k=n+p}^\infty (1-R_k)^2 \int_{\bar
 A_k} d\tilde\mu^{(2)}(\zeta)\le \frac {52}{(1-q)(1-|z|)^2}
\sum_{k=n+p}^\infty (1-R_k)= \\ =\frac {52(1-R_{n+p})}{(1-q)^2
(1-|z|)^2}\le \frac{C_2(q)}{1-R_n}.
\end{gather*}
Thus, $u_2$ is represented by the integral of the subharmonic
function $\log|E|$ of $z$, and the integral  converges uniformly
on  compact subsets in $\D$, and so $u_2\in \SH(\D)$. Since
$1-|\zeta|^2\le 3/4$ for $\zeta \in \supp \tilde \mu^{(2)}$, using
(\ref{e:pfe}) and (\ref{e:esm2}) we have
\begin{gather} \nonumber
|u_2(0)|\le  \int_{\D} |\log|E(1-|\zeta|^2, 1)||\,
d\tilde\mu^{(2)}(\zeta) \le \int_{\D} 2(1-|\zeta|^2)^2
d\tilde\mu^{(2)}(\zeta)\le
\\ \le 8\sum_{k=0}^\infty \int_{\bar A_k} (1-|\zeta|)^2 d\tilde\mu^{(2)}(\zeta)
\le
 \frac{104}{1-q} \sum_{k=0}^\infty (1-R_k)=C_3(q). \label{e:u0}
\end{gather}
Let us estimate $T(r,u_2)\stackrel{\rm def}{=} \frac 1{2\pi}
\int_0^{2\pi} u_2^+(re^{i\theta})\, d\theta$ for $r\le R_n$, where
$u^+=\max \{u, 0\}$. Note that for $|\zeta|\le R_{n+2}$, $|z|\le
R_n$ we have $\frac{1-|\zeta|^2}{|1-\bar \zeta z|}\le 2$. Thus
$$ \log \Bigr| E\Bigl( \frac{1-|\zeta|^2}{1-\bar \zeta z}, 1)\Bigl| \le
12 e \frac{1-|\zeta|^2}{|1-\bar \zeta z|}$$
in this case. Using the latter estimate, (\ref{e:pfe}), (\ref{e:esm2}), and the lemma
\cite[Ch.5.10, p.226]{Tsuji}
 we get
\begin{gather*}
T(r,u_2)\le \frac 1{2\pi}  \int_0^{2\pi} \biggl(\sum_{k=0}^{n+1} \int_{\bar A_k} 12e \frac{1-|\zeta|^2}
{|1-\bar \zeta re^{i\th}|} \, d\mu_k^{(2)}(\zeta)\biggr) d\th +\\
+\frac 1{2\pi} \int_0^{2\pi} \biggl(\sum_{k=n+2}^{\infty}\int_{\bar A_k}
6e \frac{(1-|\zeta|^2)^2}
{|1-\bar \zeta re^{i\th}|^2} \, d\mu_k^{(2)}(\zeta)\biggr) d\th \le \\
\le C_4(q) \biggl( \sum_{k=0}^{n+1}  \int_{\bar A_k} (1-|\zeta|^2) \log \frac 1{1-r}
d\mu_k^{(2)}(\zeta) +  \sum_{k=n+2}^{\infty}   \int_{\bar A_k} \frac{(1-|\zeta|^2)^2}{1-r}
d\mu_k^{(2)}(\zeta)\biggr)\le \\
\le C_5(q) \biggl( \sum_{k=0}^{n+1} \log \frac 1{1-r} + \sum_{k=n+2}^\infty  \frac {1-R_k}{1-r}
\biggr)\le
\\ \le C_6(q) n \log \frac 1{1-r} \le C_7(q) \log^2 \frac 1{1-r}.
\end{gather*}

Finally, by the First main theorem for subharmonic functions
\cite[Ch. 3.9]{Hay}
\begin{gather*}
m(r,u_2)\stackrel{\rm def}{=} \frac1{2\pi} \int_0^{2\pi} u_2^{-}
(re^{i\th}) d\th= \\ = T(r,u_2)-\int_0^r \frac {n(t,u_2)}t dt
-u_2(0)\le T(r,u_2) +C_3(q).
\end{gather*}
Therefore, $\int_0^{2\pi} |u_2
(re^{i\th})| d\th \le 4\pi T(r,u_2)+ C_8(q).$

Consequently,   $$\int_{|z|\le 1} |u_2(z)|\, dm(z)\le 4\pi
\int_0^1 T(r,u_2)\, dr+ C_8(q)\le$$ $$\le C_9(q)\int_0^1 \log^2
\frac 1{1-r} dr\le C_{10}(q).$$ Lemma 1 is proved.
\end{proof}

\subsection*{\bf 2.2. Approximation of $\tilde \mu_1$}
The following theorem plays a key role in approximation of $u$.

\begin{pheo} \label{t:ge}
 \sl Let $\mu$ be a measure in  $\R^2$ with compact support,
$\rom\supp\mu \subset \Pi$, and $\mu (\Pi)\in \N$, where $\Pi$ is
a rectangle with  ratio of  side lengths $l_0\ge 1$. Suppose, in
addition, that for any line $L$ parallel to a side of   $\Pi$,
there is at most one point $p\in L$ such that
\begin{equation} \label{e:2a}
 0<\mu(\{p\}) (<1)
\quad \mbox{while always}  \quad \mu(L\setminus \{p\})=0 ,
\end{equation}
Then there exist a system of rectangles $\Pi_k\subset \Pi$ with
sides parallel to the sides of $\Pi$, and measures $\mu_k$ with
the following properties:
\begin{itemize}
\item[1)]  $\rom\supp\mu_k \subset \Pi_k$;
\item[2)] $\mu_k(\Pi_k)=1$, $\sum_k\mu_k =\mu$;
\item[3)]  the interiors of the  convex
hulls of the supports of  $\mu_k$ are pairwise disjoint;
\item[4)] the ratio
of the side lengths of rectangles $\Pi_k$ lies  in the
interval $[1/l, l]$, where $l=\max \{l_0, 3\}$;
\item[5)] each point of the plane belongs
to the interiors of at most 4  rectangles~$\Pi_k$.
\end{itemize}
\end{pheo}
Theorem \ref{t:ge} was proved  by R.\ S.\ Yulmukhametov
\cite[Theorem 1]{Yul} for absolutely continuous measures (i.e.
$\nu$ such that $m(E)=0 \Rightarrow \nu(E)=0$) and $l_0=1$. In this case
condition (\ref{e:2a}) is fulfilled automatically. In
\cite[Theorem 2.1]{Dr} D.~Drasin showed that Yulmukhametov's proof works if
 the condition of continuity is replaced by condition
 (\ref{e:2a}). We can drop
condition (\ref{e:2a}) rotating the initial square~\cite{Dr}.
One can also  consider Theorem \ref{t:ge} as a formal consequence of
Theorem 3~\cite{AA}. Here $l_0$  plays role  for a finite set of rectangles corresponding to small $k$'s, but in \cite{AA} it plays  the principal role in the proof.

\smallskip
\noindent{\it Remark 2.1.} In the proof of Theorem \ref{t:ge}
\cite{Dr} rectangles $\Pi_k$  are obtained  by splitting the given
rectangles, starting with $\Pi$,  into smaller rectangles in the
following way. The length of the smaller side of the initial
rectangle coincides with that of a side of the  rectangle obtained
in the first generation, and the length of the other side of the
new rectangle is between one third
 and  two thirds  of the  length of the other
side of the initial rectangle. Thus we can start with a rectangle instead
of a square and $l=\max\{l_0, 3\}$.

\smallskip
Let $u_1(z)=u(z)-u_2(z)-u_3(z)$. Then $\mu_{u_1}=\tilde
\mu^{(1)}$, $\mu_{n,m}^{(1)}(\bar A_{n,m})\in 2\Z_+$, $n\in
\mathbb{Z}_+$, $0\le m\le M_n-1$.

Let  \begin{multline*} P_{n,m}=\log  {\ov A}_{n,m}=\\ =
\Bigr\{s=\si+it: \log R_n \le \si \le  \log  R_{n+1} , \frac{2\pi
m}{M_n} \le t\le \frac{2\pi (m+1)}{M_n}\Bigl\}. \end{multline*}

According to (\ref{e:22'}) the ratio of the sides of $P_{n,m}$
is
\begin{equation} \frac{\log \frac{R_{n+1}}{R_n}}{2\pi/[\frac
{2\pi}{(1-q)(1-R_n)}]} \to 1, \quad n\to\infty. \label{e:quad}
\end{equation}

Let  $d\nu_{n,m}(s)\stackrel{\rm def}=d\mu_{n,m}^{(1)}(e^s)$,
$s\in P_{n,m}$, (i.~e. $\nu_{n,m}(S)=\mu_{n,m}^{(1)}(\exp S)$ for
every Borel set $S\subset \mathbb{C}$). By  our assumptions the
conditions of Theorem \ref{t:ge} are satisfied  for
 $\Pi=P_{n,m}$ and $\mu=\nu_{n,m}/2$, and all admissible $n,m$.
By Theorem \ref{t:ge} there exists a system   ($P_{nkm}$,
$\nu_{nmk}$) of rectangles and measures,  $k\le N_{nm}$, $0\le m
\le M_n-1$ with the properties: 1)~$\nu_{nmk}(P_{nmk})=1$; 2)
$\supp \nu_{nmk} \subset P_{nmk}$; 3)~$2\sum_k \nu_{nmk}=
\nu_{n,m}$; 4)~every point $s$ such that  $\Re s<0$, $0\le  \Im s<
2\pi$ belongs to the interiors of at most four rectangles
$P_{nmk}$; 5)~the ratio of the side lengths lies between two
positive constants. Indexing the new system $(P_{nmk}, 2\nu_{nmk})$ with
the natural numbers, we
 obtain a system $(P^{(l)}, \nu^{(l)})$
with $\nu^{(l)}(P^{(l)})=2$, $\supp \nu^{(l)} \subset P^{(l)}$ etc.

Let the measure $\mu^{(l)}$ defined on $\D$ be such that
$d\mu^{(l)} (e^s)\stackrel{\rm def}=d\nu^{(l)} (s)$, $\Re s<0$,
$0\le \Im s< 2\pi$, $Q^{(l)}=\exp P^{(l)}$. Let
\begin{equation}\label{e3.2}
\zeta_l\stackrel{\rm def}=\frac12 \int\limits_{\ql} \zeta d\ml.
\end{equation}
be the center of mass of $\ql$, $l\in \N$.

We define  $\ze1$, $\ze2$ as solutions of the system
\begin{equation}
\label{e3.1}
\left\{{\aligned
\ze1+\ze2 &=\int\limits_{Q^{(l)}} \zeta d\ml, \\
(\ze1)^2+(\ze2)^2 &=\int\limits _{\ql} \zeta^2 d\ml,
\endaligned}\right.
\end{equation}
 From (\ref{e3.1}) and  (\ref{e3.2}) it follows that (see
 \cite{LM}, \cite{AA} or Lemma \ref{l:3} below)
\begin{align*}
|\ze{j}-\zeta_l|\le
\diam \ql\equiv d_l, \quad j\in\{1,2\}.
\end{align*}
Consequently
 we obtain
\begin{equation} \label{e3.3}
\max_{\zeta\in \ql} |\zeta-\ze{j}|\le 2d_l, \; j\in\{1,2\},  \quad
\sup_{\zeta\in \ql} |\zeta-\zeta_{l}|\le d_l.
\end{equation}
We write
\begin{gather} \nonumber
\Delta _l(z)\stackrel{\rm def}= \int\limits_{Q^{(l)}} \Bl \log \Bm
\frac {z-\zeta}{1-z {\bar  \zeta}} \Bm -\frac 12 \log  \Bm \frac
{z-{\ze1}} {1-z  {\bar \zeta}^{(1)}_l} \Bm - \frac 12 \log \Bm
\frac {z -{\ze2}}{1-z   {\bar\zeta^{(2)}_l}} \Bm\Br \,
d\mu^{(l)}(\zeta), \label{e3.4} \\V(z)\stackrel{\rm def}=\sum_l
\Delta _l(z). \nonumber
\end{gather}
Fix a sufficiently large $m$ (in particular, $m\ge 13$) and $z\in
A_m$.
 Let $\L^+$ be the set of $l$'s
such that $Q^{(l)}\subset U(0,R_{m-13})$, and $\L^-$ the set of
$l$'s with $Q^{(l)} \subset \{ \zeta: R_{m+13}\le |\zeta|<1\}$,
$\L^0=\N \setminus
 (\L^- \cup \L^+)$.
\begin{lemma} \label{l:zz}
There exists $l^*\in \N$ such that
 $\zeta_1$, $\zeta_2\in  U(Q^{(l)}, 2d_l)$, $l\in \L^+\cup
\L^-$, $l> l^*$ imply
$$ \frac 1{16} |z-\zeta_2|\le |z-\zeta_1|\le 16 |z-\zeta_2|.$$
\end{lemma}
\begin{proof}[Proof of Lemma \ref{l:zz}]
First, let $l\in \L^+$, i.e. $z\in A_m$, $\ql\subset \ov{A_p}$,
$p\le m-13$. In view of (\ref{e:quad}) $\ql=\exp\pl $ is ``almost
a square''. More precisely, there exists $l^* \in \N$ such that
for all $l>l^*$ $$ \diam \ql=d_l< \frac32 (R_{p+1} -R_p), \quad
\ql \subset \ov{A_p}.$$ Since  $\zeta_1$, $\zeta_2\in  U(Q^{(l)},
2d_l)$, we have
\begin{gather}
\label{e:zz0} R_p-3(R_{p+1}-R_p)\le |\zeta_2|\le
R_{p+1}+3(R_{p+1}-R_p),
 \\ \label{e:zz1} |z-\zeta_1|\ge
|z-\zeta_2|-|\zeta_2-\zeta_1|\ge |z-\zeta_2| -5 d_l \ge
|z-\zeta_2|- \frac {15}2 (R_{p+1}-R_p).
\end{gather}
On the other hand, by the choice of $q$ (see (\ref{e:qcon})) and
(\ref{e:zz0})
\begin{gather*}
|z-\zeta_2| \ge R_m-R_{p+1}\ge R_{p+13}-R_{p+1}-3(R_{p+1}- R_p)=\\
=\sum_{s=1}^{12} (R_{p+s+1} -R_{p+s})-3(R_{p+1}- R_p)=
\Bigl(\sum_{s=1}^{12} q^{s}-3\Bigr)(R_{p+1}- R_p)>8 (R_{p+1}-R_p).
\end{gather*}
The latter inequality and (\ref{e:zz1}) yield
\begin{gather*}
 |z-\zeta_1|\ge |z-\zeta_2| -\frac{15}2
(R_{p+1}-R_p)> |z-\zeta_2| -\frac {15}{16} |z-\zeta_2|=\frac
1{16}|z-\zeta_2|.
\end{gather*}

For $l\in \L^-$, $\ql\subset \{R_{m+13}\le |\zeta|<1\}$ we have
$p\ge m+13$, and
inequality (\ref{e:zz1}) still holds.

Similarly, by the choice of $q$ and (\ref{e:zz0})
\begin{gather*}
|z-\zeta_2| \ge R_{p-3} -R_{m+1}\ge R_{p-3}-R_{p-12}
\ge 9q^4 (R_{p+1}- R_p)>8 (R_{p+1}-R_p),
\end{gather*}
that together with (\ref{e:zz1}) implies required inequality in
this case. Lemma \ref{l:zz} is proved.
\end{proof}

Let $l\in \L^-\cup \L^+$.
 For
$\zeta\in \ql$, we define $L(\zeta)=L_l(\zeta)=\log  \bigl(
\frac{z-\zeta}{1-\bar z \zeta}\bigr)$, where $\log w$ is an
arbitrary branch of  $\mathop{\rm Log} w$ in $w(\ql)$, $w(\zeta)=\frac{z-\zeta}{1-\bar z \zeta}$. Then
$L(\zeta)$ is analytic in $\ql$.
 We shall use the following identities
\begin{align}\nonumber
&L(\zeta)-L(\ze1)=\int\limits_{\ze1}^\zeta L'(s)\, ds =
L'(\ze1)(\zeta-\ze1) + \int\limits_{\ze1}^\zeta L''(s)(\zeta-s)\,
ds= \nonumber
\\
&=L'(\ze1)(\zeta-\ze1)
+\frac12 L''(\ze1)(\zeta-\ze1)^2
+\frac 12\int\limits_{\ze1}^\zeta L'''(s)(\zeta-s)^2\, ds.
\label{e3.5}
\end{align}

Elementary geometric arguments show that $|\frac 1{\bar z}
-\zeta|^{-1} \le |z-\zeta|^{-1}$ for $z,\zeta\in \D$.  Since
$L'(\zeta)=\frac{1}{\zeta -z} +\frac {\bar z}{1-\bar z \zeta}$, we
have \bey \label{e:lde} |L'(\zeta)|\le \frac 2{|\zeta-z|}, \quad
|L''(\zeta)|\le \frac 2{|\zeta-z|^2}, \quad |L'''(\zeta)|\le \frac
4{|\zeta-z|^3}. \eey

Now  we estimate
 $|\De_l(z)| $ for $l\in \L^+ \cup \L^-$. By the definitions of
$L(\zeta)$,  $\De_l(z)$, (\ref{e3.5}) and (\ref{e3.1}) 
 we have
\begin{gather} \nonumber
|\Delta_l(z)|=\Bm \Re \int\limits_{\ql} \Bl L(\zeta
)-L(\ze1)-\frac12 (L(\ze2)-L(\ze1)\Br \, d\ml\Bm = \\ =\Bm \Re
\int\limits_{\ql} \biggl( L'(\ze1)\Bl \zeta -
\frac12(\ze1+\ze2)\Br + \nonumber\\ +\int_{\ze1}^\zeta
L''(s)(\zeta-s)ds -\frac 12 \int_{\ze1}^{\ze2} L''(s)(\ze2-s)ds
\biggr) \, d\ml\Bm=\nonumber \\ =\Bm \Re \int\limits_{\ql}\biggl(
\int_{\ze1}^\zeta L''(s)(\zeta-s)ds -\frac 12 \int_{\ze1}^{\ze2}
L''(s)(\ze2-s)ds \biggr) \, d\ml\Bm  \label{e:2.14'}\end{gather}
Using (\ref{e:2.14'}), (\ref{e:lde}) and (\ref{e3.3}), we obtain
\begin{gather} \nonumber
|\Delta_l(z)|\le \int\limits_{\ql} \int\limits_{\ze1}^\zeta \frac
{2|\zeta-s|}{|s-z|^2} |ds|\, d\ml+\frac 12 \int\limits_{\ql}
\int\limits_ {\ze1}^{\ze2} \frac {2|\ze2 -s||ds| }{|s-z|^2} \,
d\ml\le \nonumber
\\
\le  {12 d_l^2}\max_{s\in B_l}  \frac1{|s-z|^2}, \label{e3.6}
\end{gather}
where $B_l=\ov{U(\ql, 2d_l)}$. Applying Lemma \ref{l:zz}, we have
$(z\in \bar A_m)$
\begin{gather} \nonumber
\sum_{l\in \L^-} |\De_l(z)| \le 12  \sum_{l\in {\L}^-}  d_l^2
\max_{s\in B_l} \frac1{|s-z|^2}\le  C_{11} \sum_{l\in {\L}^-}
\int_{\ql} \frac {dm(z)}{|z-\zeta|^2} \le \\ \le 4C_{11}
\int_{R_{m+13}\le |\zeta |<1}  \frac {dm(z)}{|z-\zeta|^2} \le
C_{12} \int_{R_{m+13}}^1 \frac {d\rho}{\rho-|z|} \le C_{13}(q).
\label{e:smi}
\end{gather}
Similarly,
\begin{gather} \nonumber
\sum_{l\in \L^+} |\De_l(z)| \le 12  \sum_{l\in {\L}^+}  d_l^2
\max_{s\in B_l} \frac1{|s-z|^2}\le 4 C_{11} \sum_{l\in {\L}^+}
 \int_{|\zeta|\le R_{m-13} }  \frac
{dm(z)}{|z-\zeta|^2} \le \\ \le C_{12} \int_0^{R_{m-13}} \frac
{d\rho}{|z|-\rho} \le C_{14}(q) \log \frac 1{1-|z|}. \label{e:spl}
\end{gather}
Hence, \bey \label{e:spm} \int_{|z|\le R_n} \sum_{l\in \L^+\cup
\L^-} |\De_l(z)|\, dm(z) < C_{15}(q). \eey

It remains to estimate $\int_{|z|\le R_n} \sum _{l\in
\L^0}|\De_l(z)|\, dm(z) $. Here we follow the arguments from
\cite[e.-g.]{LM}. If $\dist (z, \ql)>10 d_l$, similarly to
(\ref{e:2.14'}), from  (\ref{e3.5}), (\ref{e3.1}), (\ref{e:lde}) and
(\ref{e3.3}) we deduce
\begin{gather} \nonumber
|\De_l(z)| = \Bm \Re \int\limits_{\ql} \biggl( L'(\ze1)\Bl \zeta -
\frac12(\ze1+\ze2)\Br + \nonumber\\ + \frac {
L''(\ze1)}2\Bl\zeta^2 -\frac{(\ze1)^2+(\ze2)^2}2
+\ze1(\ze1+\ze2-2\zeta)\Br + \nonumber \\ +
\frac12\int_{\ze1}^\zeta L'''(s)(\zeta-s)^2ds  -\frac 14
\int_{\ze1}^{\ze2} L'''(s)(\ze2-s)^2ds \biggr) \,
d\ml\Bm=\nonumber
\\=\biggl|  \Re \int\limits_{\ql} \biggl( \frac 12
\int\limits_{\ze1}^{\zeta} L'''(s)(\zeta-s)^2\, ds -\frac14
\int\limits_{\ze1}^{\ze2} L'''(s)(\zeta-s)^2\, ds \biggr)\, d\ml
\biggl|\le \nonumber
\\
\le  6 d_l^3  \max_{s\in B_l}  \frac{1}{|s-z|^3} \le \frac
{6d_l^3}{| \zeta_l^{(1)}-z|^3}\max _{s\in B_l} \Bl 1+\frac
{|\ze1-s|}{|s-z|}\Br^3\le \frac {26d_l^3}{| \zeta_l^{(1)}-z|^3}.
\label{e:dlfar}
\end{gather}

Since $\L_0$ depends only on $m$ when $z\in A_m$, we have
\begin{gather} \nonumber
\int\limits_{\ov{A}_m} \sum_{l\in \L^0} |\De_l(z)| \, dm(z)\le
\sum_{l\in \L^0} \biggl( \int\limits_{\ov{A}_m \setminus U( \ze1,
10d_l)} +\int\limits_{U(\ze1, 10d_l)} \biggr) |\De_l(z)|\, dm(z)
\le \\ \le \sum_{l\in \L^0} \biggl( \int\limits_{\ov{A}_m
\setminus U(\ze1, 10d_l)} \frac{26d_l^3}{|z-\ze1|^3} dm(z)
+\int\limits_{U(\ze1, 10d_l)} |\De_l(z)|\, dm(z) \biggr).
\label{e:2.19}
\end{gather}
For the first sum we obtain
\begin{gather} \nonumber
\sum_{l\in \L^0} 26d_l^3 \int_{\ov{A}_m \setminus U(\ze1, 10d_l)}
\frac{1}{|z- \ze1|^3} dm(z) \le \sum_{l\in\L^0} 52\pi d_l^3
\int_{10d_l}^2 \frac {tdt}{t^3} \le
\\ \le 6\pi \sum_{l\in \L^0} d_l^2 \le C_{16} \sum_{l\in \L^0}
m(\ql). \label{e:2.20}
\end{gather}
We now estimate the second sum. By the definition of $\De_l(z)$
\begin{gather*}
\De_l(z)=\int_{\ql} \Bl \log \Bm \frac{z-\zeta}{10d_l} \Bm -\frac
12 \log \Bm \frac{z-\ze{1}}{10d_l} \Bm -
\frac 12 \log \Bm \frac{z-\ze{2}}{10d_l} \Bm
\Br \, d\ml- \\
-\int_{\ql} \Bigl(\log|1-z\ov{\zeta}|-\frac 12 \log|1-z\ov{\ze{1}}|
-\frac12 \log|1-z\ov{\ze{2}}|\Bigr)\, d\ml\equiv I_1+I_2.
\end{gather*}
The integral $\int |I_1| \, dm(z)$ is estimated  in \cite[g.]{LM},
\cite[p.232]{AA}. We have \bey\label{e:2.21} \int_{U(\ze1, 10d_l)}
|I_1| \, dm(z) \le C_{17} m(\ql). \eey

To estimate $|I_2|$ we note that for $l$ sufficiently large,
$|z-\zeta|\le 15 d_l$, $ \zeta\in U(\ql, 2d_l)$, $z\in \D$  we
have $|\arg z- \arg \zeta|\le 16 d_l\le 16(1-|z|)|$ by the choice of $q$.
Hence,
$$|\frac1z-\bar \zeta|\le \frac 1{|z|} -1 +1 -|\zeta|+ |\zeta||1-e^{i(\arg \zeta-\arg z)}|\le
C_{17}'(1-|z|).$$
Thus, $|1/z-\bar \zeta|\asymp 1-|z|$.
Therefore $$|I_2|\le \int_{\ql} \frac 12 \Bm \log  \frac{|\frac 1z
-\bar\zeta|^2}{|\frac 1z-\ze1||\frac 1z-\ze2|} \Bm d\ml \le
C_{18}. $$ Thus, \bey \label{e:2.22} \int_{U(\ze1,
10d_l)} |I_2| \, dm(z) \le C_{19}(q) m(\ql). \eey Finnaly, using
(\ref{e:2.20})--(\ref{e:2.22}) we deduce
\begin{gather*}
\int_{\bar A_m} \sum_{l\in \L^0} |\De_l(z)| \, dm(z) \le \\ \le
C_{20} \sum_{l\in \L^0}  m(\ql) \le 4\pi C_{20} (R_{m+13}^2
-R_{m-13}^2)\le C_{21}(q) (R_{m+1} -R_m) .
\end{gather*}

Hence, $\int_{|z|\le R_n} \sum_{l\in \L^0} |\De_l(z)|dm(z)\le
C_{20}(q)$, and this with  (\ref{e:spm}) yields that
\bey\label{e:2.23} \int_{|z|\le R_n}  |V(z)|dm(z)\le C_{22}(q),
n\to+\infty. \eey

Now we construct the function $f_1$ approximating $u_1$.

Let $K_n(z)=u_1(z)-\sum_{\ql\subset \ov{U(0,R_n)}} \De_l(z),$
$K(z)=u_1(z)-V(z)$. By the definition of $\De_l(z)$, $K_n\in
\SH(\D)$ and $$\mu_{K_n} \bigr|_{U(0,R_n)} (z)= \sum_{l=1}^n
\bigl(\delta(z-\ze1)+\de (z-\ze2)\bigr)$$ where $\delta(\zeta)$ is the unit mass supported at $u=0$. For $|z|\le R_n$, $j\ge
N\ge n+14$ as in (\ref{e:smi}) we have
\begin{gather*}
|K_j(z)-K(z)| \le \sum_{\ql\subset \{|\zeta|\ge R_{N+1}\}}
|\De_l(z)| \le \\ \le  C_{23} \int_{R_{N+1}\le |\zeta| <1} \frac
{dm(z)}{|z-\zeta|^2} \le C_{24} \frac {1-R_{N+1}}{R_{N+1} -|z|}
\to 0, \quad N\to +\infty.
\end{gather*}
Therefore $K_n(z)\rightrightarrows K(z)$ on the compact sets in
$\D$ as $n\to+\infty$, and $\mu_K\Bm_{\D}
=\sum_{l}(\delta(z-\ze1)+\de (z-\ze2)$. Hence,
$K(z)=\log|f_1(z)|$, where $f_1$ is analytic in $\D$.

\subsection*{2.3. Approximation of $u_3$}

Let $u_3$ be defined by (\ref{e:u3}), $$N=2\bigl[n(1/2,u_3)/2
\bigr], \quad \rho_0=\inf\{ r\ge 0: n(r, u_3)\ge N\}.$$ We
represent $\mu_{u_3} =\mu^1+ \mu^2$ where $\mu^1$ and $\mu^2$ are
measures such that
\begin{gather*}
\supp \mu^1\subset \ov{U(0,\rho_0)}, \quad \supp \mu^2\subset
\ov{U\Bigl(0, \frac12\Bigr)} \setminus U(0,\rho_0), \\
\mu^1\Bigl({U\Bigl(0,\frac 12\Bigr)\Bigr)}=N, \quad 0\le
\mu^2\Bigl({U\Bigl(0,\frac 12\Bigl)\Bigl)}<2.
\end{gather*}

Let $v_2(z)=\int\limits_{{U(0, \frac 12)}} \log |z-\zeta|\,
d\mu^2(\zeta)$. Then, using the last estimate,
\begin{gather*}
\int\limits_{\D} | v_2(z)|\, dm(z)\le
\int\limits_{{U(0,1/2)}}\int\limits_{\D} |\log|z-\zeta||\, dm(z)\,
d\mu^2(\zeta) \le \\ \le \int\limits_{{U(0,1/2)}}
\int\limits_{U(\zeta,2)} |\log|z-\zeta||\, dm(z)\, d\mu^2(\zeta)
\le C_{25} n\Bigl(\frac 12, v_2\Bigr)\le 2C_{25}.
\end{gather*}

If $N=0$ there remains nothing to prove. Otherwise, we have to
approximate
\begin{equation}\label{e:v1}
  v_1(z)=u_3(z)-v_2(z)=\int\limits_{\ov{U(0,\rho_0)}} \log |z-\zeta|\,
  d\mu^1(\zeta).
\end{equation}
In this connection we recall the question of Sodin  (Question 2 in
\cite[p.315]{Sod}).

Given a Borel measure $\mu$ we define the  logarithmic potential
of~$\mu$ by the equality $$\mathcal{U}_\mu(z)=\int\log|z-\zeta|\,
d\mu(\zeta).$$

\noindent{\bf Question.} {\it Let $\mu$ be a probability measure
supported by the square $\mathcal{Q}=\{z=x+iy: |x|\le \frac 12,
|y|\le \frac 12\}$. Is it possible to find a sequence of
polynomials $\mathcal{P}_n$, $\deg \mathcal{P}_n=n$, such that $$
\iint \limits_{\begin{substack} {{|x|\le 1}
\\ {|y|\le 1} }
\end{substack}} |n\mathcal{U}_\mu(z) -\log |\mathcal{P}_n(z)||\, dxdy=O(1) \; (n\to+\infty)? $$
}

We should say that the solution is given essentially in \cite{LM}, but not asserted. Hence we prove the following

\noindent {\bf Proposition.}
 {\it Let $\mu$ be a measure supported by the square $\mathcal{Q}$, and $\mu(\mathcal{Q})=N\in \mathbb{N}$.
 Then there is an absolute constant $C$ and a polynomial $P_N$ such that
 $$ \iint \limits_{\Xi} |\mathcal{U}_\mu(z) -\log |\mathcal{P}_N(z)||\, dxdy< C, $$
 where $\Xi=\{ z=x+iy: {|x|\le 1},  {|y|\le 1}\}$.}

\begin{proof}[Proof of the proposition]
As in the proof of Theorem \ref{t:1}, if there are points $p\in
\mathcal{Q}$ such that $\mu(\{p\})\ge 1$ we represent  $\mu=\nu+
\tilde \nu$ where for any $p\in \mathcal{Q}$ we have
$\nu(\{p\})<1$, and $\tilde \nu$ is a finite (at most $N$ summand)
sum of the Dirac measures. Then $\mathcal{U}_{\tilde \nu}
=\log\prod_k |z-p_k|$, so it remains to approximate $\mathcal{U}_\nu$. By
Lemma 2.4 \cite{Dr} there exists a rotation to the system of
orthogonal coordinates such that if $L$ is any line parallel to
either of the coordinate axes, there is at most one point $p\in L$
with $\nu(\{p\})>0$, while always $\nu(L\setminus \{p\})=0$.
After the rotation the support of the new measure, which is still denoted by $\nu$,
 is contained in
$\sqrt{2} \mathcal{Q}$.

If $\omega$ is a probability measure supported on $\mathcal{Q}$,
then $\iint_\Xi |\mathcal{U}_\omega(z)| dm(z)$ is uniformly
bounded. Therefore we can assume that $N\in 2\mathbb{N}$.

By Theorem E there exists a system $(P_l, \nu_l)$ of rectangles
and measures $1\le l\le M_\nu$ with the properties:
1)~$\nu_l(P_l)=2$; 2)~$\supp \nu_l\subset P_l$; 3)~$\sum _l \nu_l
=\nu$; 4)~every point $s\in \mathcal{Q}$ belongs to interiors of at most
four rectangles $P_l$; 5)~ratio of  side lengths lays between 1/3
and 3.

Let
\begin{equation}\label{e:4.1}
\xi_l=\frac12 \int\limits_{P_l} \xi d\nu_l(\xi).
\end{equation}
be the center of mass of $P_l$, $1\le l\le M_\nu$.

We define  $\xxi1$, $\xxi2$ as solutions of the system
\begin{equation*}
\left\{{\aligned \xxi1+\xxi2 &=\int\limits_{P_{l}} \xi
d\nu_l(\xi),
\\ (\xxi1)^2+(\xxi2)^2 &=\int\limits _{P_l} \xi^2 d\nu_l(\xi),
\endaligned}\right.
\end{equation*}
We have
\begin{align*}
|\xxi{j}-\xi_l|\le \diam P_l\equiv D_l, \quad j\in\{1,2\},
\end{align*}
\begin{equation} \label{e:4.3}
\max_{\xi\in P_l} |\xi-\xxi{j}|\le 2D_l, \; j\in\{1,2\}, \quad
\sup_{\xi\in P_l} |\xi-\xi_{l}|\le D_l.
\end{equation}
We write
\begin{gather} \nonumber
\Omega(z)=\sum_l \int\limits_{P_{l}} \Bl \log  \Bm  {z-\xi} \Bm
-\frac 12 \log  \Bm {z-{\xxi1}}  \Bm - \frac 12 \log  \Bm  {z
-{\xxi2}} \Bm\Br \, d\nu_{l}(\xi) \equiv \\ \equiv \sum_l \delta
_l(z). \label{e:4.4}
\end{gather}

Since we have rotated the system of coordinate, it is sufficient
to prove that $\int_{\ov{U(0,\sqrt{2})}} |\Omega(z)|\, dm(z)$ is
bounded by an absolute constant.

For $\xi \in P_l$, $z\not\in P_l$ we define
$\lambda(\xi)=\lambda_l(\xi)=\log \bigl( {z-\xi}\bigr)$, where
$\log (z-\xi)$ is an arbitrary branch of $\mathop{\rm Log}
(z-\xi)$  in $z-P_l$. Then $\lambda(\xi)$ is analytic in $P_l$.

We have \bey \label{e:lde1}
 |\lambda'''(\xi)|\le
\frac 2{|\xi-z|^3}. \eey

As in subsection 2.2 we have
\begin{gather}\nonumber
|\delta_l(z)|\le \Bm \Re \int\limits_{P_l} \Bl \lambda(\xi
)-\lambda(\xxi1)-\frac12 (\lambda(\xxi2)-\lambda(\xxi1)\Br \,
d\nu_l(\xi)\Bm \le \\ \le \biggl|  \Re \int\limits_{P_l} \biggl(
\frac 12 \int\limits_{\xxi1}^{\xi} \lambda'''(s)(\xi-s)^2\, ds
-\frac14 \int\limits_{\xxi1}^{\xxi2} \lambda'''(s)(\xi-s)^2\, ds
\biggr)\, d\nu_l(\xi) \biggl|. \label{e:4.20}
\end{gather}
 If $\dist (z, P_l)>10 D_l$ the last estimate and (\ref{e:lde1})
 yield
\begin{gather} \nonumber
|\de_l(z)| \le  24 D_l^3  \max_{s\in E_l}  \frac{1}{|s-z|^3} \le
\frac {24D_l^3}{|\xi_l^{(1)}-z|^3}\max _{s\in E_l} \Bl
1+\frac{|\xi_l^{(1)} -s|}{|s-z|}\Br\le\frac {103
D_l^3}{|\xi_l^{(1)}-z|^3}, \label{e:dlfar1}
\end{gather}
where $E_l=\ov{U(P_l, 2D_l)}$.

Then
\begin{gather*} \nonumber
 \int\limits_{\ov{U(0,\sqrt{2})} \setminus  U(\xi_l^{(1)}, 10D_l)}
\frac{103 D_l^3}{|z- \xi_l^{(1)}|^3} dm(z) \le 206 \pi D_l^3
\int\limits_{10D_l}^2 \frac {tdt}{t^3} \le \\ \le 21\pi D_l^2 \le
C_{26} m(P_l).
\end{gather*}
On the other hand, by the definition of $\de_l(z)$
\begin{gather*}
\int\limits_{U( \xi_l^{(1)}, 10D_l)}
\de_l(z)dm(z)=\int\limits_{U(\xi_l^{(1)}, 10D_l)} \int_{P_l} \Bl
\log \Bm \frac{z-\xi}{10D_l} \Bm -\\ -\frac 12 \log \Bm
\frac{z-\xxi{1}}{10D_l} \Bm - \frac 12 \log \Bm
\frac{z-\xxi{2}}{10D_l} \Bm \Br \, d\nu_l(\xi) dm(z)\le C_{27}
m(P_l).
\end{gather*}
 From (\ref{e:4.4}) and  the latter estimates, it follows that
\begin{equation}\label{e:4last}
\int\limits_{\ov{U(0,\sqrt{2})}}  |\Omega(z)| dm(z)\le   \sum_{l} \int\limits_{\ov{U(0,\sqrt{2})}} \delta_l(z)
dm(z)\le C_{28} \sum_l m(P_l)\le 4 C_{28}m(\sqrt{2}\mathcal{Q})=C_{29}.
\end{equation}
Thus, $\mathcal{P}(z)=\prod_l (z-\xi_l^{(1)})(z-\xi_l^{(2)})$ is a
required  polynomial. This completes the proof of the proposition.
\end{proof}

Finally, let $f=f_1 \mathcal{P}$. By Lemma \ref{l:d1},
(\ref{e:2.23}), and (\ref{e:4last}) we have $(n\to+\infty)$
\begin{gather*}
\int_{|z|\le R_n} |\log|f(z)|-u(z)||\, dm(z) \le \int_{|z|\le R_n}
(|K(z)-u_1(z)|| + |u_2(z)| +\\ + |\log |\mathcal{P}|-u_3(z)|)\,
dm(z) \le \int_{|z|\le R_n} (|V(z)|+|\Omega(z)|)\, dm(z)+
C_{10}(q)\le C_{30}(q).
\end{gather*}
Fixing any $q$ satisfying (\ref{e:qcon}) we finish the proof of
Theorem \ref{t:1}.

\section{Uniform approximation} \label{s:3}
In this section we prove some counterparts of results due to
Yu.Lyubarskii and Eu.Malinnikova \cite{LM}. We start with
counterparts of notions introduced in \cite{LM}, which reflect
regularity properties of measures.

\noindent {\bf Definition 1.} {Let $b\colon [0,1)\to(0,+\infty)$
be such that $b(r)\le 1-r$, \bey b(r_1)\asymp b(r_2) \quad as
\quad  1-r_1\asymp 1-r_2,\quad r_1\uparrow 1. \label{e:sgb} \eey A
measure $\mu$ on $\mathbb{D}$ admits a {\it partition of slow
variation with the function $b$} if there exist integers $N$, $p$
and sequences $(\ql)$ of subsets of $\D$ and $(\mu^{(l)})$ of
measures with the following properties:
\begin{itemize}
\item[i)] $\supp \mu^{(l)} \subset \ql$, $\mu^{(l)}(\ql)=p$;
\item[ii)] $\supp (\mu -\sum_l \mu^{(l)})  \subset \D$, $(\mu-\sum_l
\mu^{(l)})(\D)<+\infty$.
\item[iii)] $1-\dist (0, \ql)\ge K(p)\diam \ql$, and each $z\in \D$ belongs to at most $N$ various
$\ql$'s;
\item[iv)] For each  $l$ the set $\log \ql$ is  a rectangle with
sides parallel to the coordinate axes, and the ratio of
sides lengths lies between two  positive constants independent of
$l$.
\item[v)] $\diam \ql \asymp b(\dist (\ql, 0))$.
\end{itemize}
}

\noindent{\it Remark 3.1.}  This is similar to  \cite{LM}, except
we have introduced the parameter $p$ ($p=2$ in \cite{LM}). Property iii) is
adapted for $\D$.

\noindent{\bf Definition 2.} {Given a function $b$ satisfying
(\ref{e:sgb}) we say that a {\it  measure $\mu$ is locally regular
with respect to (w.r.t.) $b$ if $$ \int_0^{b(|z|)} \frac
{\mu(U(z,t))}t dt =O(1), \quad r_0<|z|<1, $$ for some constant
$r_0\in (0,1)$.}

\begin{theorem} \label{t:3} \it
Let $u\in \SH(D)$, $b\colon [0,1)\to (0,+\infty)$ satisfy
(\ref{e:sgb}). Let $\mu_u$ admits a partition of slow variation, assume that
$\mu_u$ is locally regular w.r.t. $b$, and, with $p$ from above, that \bey \label{e:bco}
\int_0^1 \frac {b^{p-1}(t)}{(1-t)^p} dt<+\infty. \eey Then there
exists an analytic function $f$ in $\D$ such that
$\forall\varepsilon>0$ $\exists r_1\in (0,1)$ $$ \log|f(z)|
-u(z)=O(1), \quad r_1<|z|<1, \; z\not\in E_\ve $$ where $E_{\ve}
=\{ z\in \D: \dist (z, Z_f)\le \ve b(|z|)\}$, and for some
constant $C>0$ \bey \label{e:bo_ab} \log|f(z)| -u(z)< C, \quad
z\in \D. \eey Moreover, $Z_f \subset U(\supp \mu_u, K_1(p)
b(|z|))$, $K_1(p)$ is a positive constant, and
\begin{equation}\label{e:tt}
  T(r,u)-T(r,\log |f|)=O(1), \quad r\uparrow1.
\end{equation}
\end{theorem}

\noindent{\it Remark 3.2.} The author does not know whether
condition (\ref{e:bco}) is necessary. But if
$b(t)=O((1-t)\log^{-\eta} (1-t)$, $\eta>0$, $t\uparrow1$
(\ref{e:bco}) holds for sufficiently large $p$. On the other hand,
in view of v) the condition $b(t)=O(1-t)$ as $t\uparrow 1$ is
natural.

\begin{proof}[Proof of Theorem \ref{t:3}]
We follow \cite{LM} and also use arguments and notation from the
proof of Theorem \ref{t:1}.

Let $\tilde \mu =\mu_u - \sum_{l} \mu^{(l)}$. Since $\Bm \frac
{z-\zeta}{1-\bar \zeta z}\Bm \to 1$ as $ |z|\uparrow1$ for fixed
$\zeta \in \D$, $\tilde \mu (\D)<+\infty$, $$ \tilde
u_1(z)=\int_{\D} \log \Bm \frac {z-\zeta}{1-\bar \zeta z}\Bm\,
d\tilde \mu(\zeta)$$ is a subharmonic function in $\D$ and
$|\tilde u_1(z)|<C$ for $r_1<|z|<1$, $r_1\in (0,1)$. So we can
assume that $\mu_u=\sum_l \mu^{(l)}$ where $\mu^{(l)}$ are from
Definition 1.

Fix a partition of slow variation. Instead of points $\ze1$ and
$\ze2$ satisfying (\ref{e3.1}) we define $\xi_1^{(l)}$, \dots, $\xi_{p}^{(l)}$
from the system
\begin{equation} \label{e:xisys}
\begin{cases} \xi_1+\dots+\xi_p&=\int_{\ql} \xi d\mu^{(l)}(\xi), \\
\xi_1^2+\dots+\xi_p^2&=\int_{\ql} \xi^2 d\mu^{(l)}(\xi), \\ &
\vdots
\\ \xi_1^p+\dots+\xi_p^p&=\int_{\ql} \xi^p d\mu^{(l)}(\xi),
\end{cases}
\end{equation}
Lemma \ref{l:3} is a modification of the estimates in(\ref{e3.3}).
\begin{lemma} \label{l:3}
Let $\Pi$ be a set in $\C$, $\mu$ is a measure on $\Pi$,
$\mu(\Pi)=p\in \N$, $\diam \Pi=d$. Then for any solution $(\xi_1,
\dots, \xi_p)$ of (\ref{e:xisys}) we have $|\xi_j-\xi_0|\le K_1(p)
d$ where $ K_1(p)$ is a constant, $\xi_0$ is the center of mass of
$\Pi$.
\end{lemma}

\begin{proof}[Proof of Lemma \ref{l:3}]
Let $\xi_0=\frac 1p \int_\Pi \xi d\mu(\xi)$ be the center of mass
of $\Pi$.
By  induction, it is easy to prove that  (\ref{e:xisys}) is
equivalent  to the system
\begin{equation} \label{e:wsys}
\begin{cases} w_1+&\dots+w_p=0, \\
w_1^2+&\dots+w_p^2=J_2, \\ &\vdots \\ w_1^p+&\dots+w_p^p=J_p,
\end{cases}
\end{equation}
where $w_k=\xi_k-\xi_0$, $J_k=\int_\Pi (\xi-\xi_0)^k\, d\mu(\xi)$,
$1\le k\le p$. Note that $$|J_k|\le \int_\Pi |\xi-\xi_0|^k\,
d\mu(\xi) \le p d^k.$$ From algebra it is well-known that the
symmetric polynomials $$\sum_{1\le i_1< \dots <i_k\le m} w_{i_1}
\cdots w_{i_k},$$  $1\le k \le m$, can be obtained  from the
polynomials $\sum_{j=1}^m w_j^k$ using only finite number of
operations of addition and multiplication. Therefore
(\ref{e:wsys}) yields
\begin{equation*} 
\begin{cases} w_1+\dots+w_p=0, \\
\sum\limits_{1\le i_1< i_2\le p} w_{i_1} w_{i_2}=b_2 , \\
\hphantom{ssss}\vdots
\\ w_1 \cdots w_p=b_p,
\end{cases}
\end{equation*}
where $b_k=\sum_l a_{lk} (J_1)^{s_{1l}^{(k)}} \cdots
(J_m)^{s_{ml}^{(k)}}$,  $a_{lk} =a_{lk}(p)$, $s_{jl}^{(k)}$ are
non-negative integers, and $\sum_{j=1}^p s_{jl}^{(k)}j=k$. The
last equality follows from homogeneousity. Hence, there exists a
constant $K_1(p)\ge 2$ such that $|b_k|\le K_1(p) d^k$, $1\le k\le
p$. By  Vieta's  formulas (\cite[\S\S 51,52]{Ku}) $w_j$, $1\le j\le p$, satisfy the
equation
\begin{equation} \label{e:weq}
w^p +b_2w^{p-2} -b_3w^{p-3} +\dots +(-1)^p b_p=0.
\end{equation}
For $|w|=K_1(p) d$ we have
\begin{gather*}
|w^p +b_2w^{p-2} -b_3w^{p-3} +\dots +(-1)^p b_p|\le K_1(p)(d^2
|w|^{p-2} +\dots +d^p)= \\ = K_1(p) d^p (K_1^{p-2} + K_1^{p-1}
+\dots +1) < 2 K_1^{p-1}(p) d^p\le  K_1^p(p) d^p =|w|^p.
\end{gather*}
By  Rouch\'e's theorem all $p$ roots of  (\ref{e:weq}) lay in the
disk $|w|\le K_1(p) d$, i.e. $|\xi_j -\xi_0|\le K_1(p) d$.
Consequently, $\dist (\xi_j, \Pi)\le K_1(p)d$.
\end{proof}
Applying  Lemma \ref{l:3} to  $\ql$ we obtain that
$|\xi_l^{(j)} -\xi_l|\le K_1(p) d_l$, $1\le j\le p$, where
$\xi_l=\frac 1p \int_{\ql} \xi d\mu^{(l)}(\xi).$

Consider $$ V(z)=\sum_l j_l(z)\stackrel{\rm def}= \sum_{l}
\int_{\ql} \Bl \log \Bm \frac{z-\zeta}{1-\bar z\zeta}\Bm -\frac 1p
\sum_{j=1}^p \log \Bm\frac{z-\xi_l^{(j)}}{1-\bar
z\xi_l^{(j)}}\Bm\Br d\mu^{(l)} (\zeta). $$ For $R_n=1-2^{-n}$, $z
\in A_m$, $m$ is fixed, we define sets of indices $\L^+$, $\L^-$
and $\L^0$ as in the proof of Theorem \ref{t:1}.

The estimate of $\sum\limits_{l\in \L^-} j_l(z)$ repeats that of
$\sum\limits_{l\in \L^-} \De_l(z)$, so \bey \sum_{l\in \L^-}
|j_l(z)| \le C_{31}. \label{e:ejm} \eey Following \cite{LM} we
estimate $\sum_{l\in \L^0} j_l(z)$. Let $b_m=b(R_m)$. Note that
$d_l\asymp b_m $ for $l\in \L^0$ by condition v). As in
(\ref{e3.6}) we have
\bey \label{e:ejl0} |j_l(z)|\le C_{32}d_l^3
\max _{s\in \ov{U(\ql, K_1(p)d_l)}} |s-z|^{-3}\le C_{32}'\frac
{d_l^3}{|\xi_l^{(1)} -z|^3}, \eey
 provided that $\dist (z,\ql)\ge
3K_1(p)d_l$. Then
\begin{gather} \nonumber
\biggl| \sum_{\begin{substack} {l\in \L^0 \\ \ql \cap U(z,3K_1(p)d_l)=\varnothing}
\end{substack}} j_l(z) \biggr| \le C_{32} \sum_{l\in \L^0} \frac
{d_l^3}{| \xi_l^{(1)} -z|^3}  \le  \\ \le C_{33} b_m
\int\limits_{|z-\zeta|>C_{34}b_m} \frac {dm(\zeta)} {|z-\zeta|^3}
\le C_{35} \frac {b_m}{b_m} =C_{35}. \label{e:ej0}
\end{gather}
Let now $l$ be such that
$\ql \cap U(z, 3K_1(p)d_l)\ne \varnothing$. Since $d_l\asymp b_m$, the
number of these $l$ is bounded uniformly in $l$.
For $z\not \in E_\ve$ we have $(1\le k\le p)$
\bey \label{e:zmz}
\log|z-\xi_l^{(k)}| =\log b_m + \log \frac{|z-\xi_l^{(k)}|}{b_m} =\log
b(|z|)+O(1).
\eey
Therefore
\begin{gather*}
j_l(z)=\int_{\ql} \Bl \log|z-\zeta| -\frac 1p\sum_{k=1}^p
\log|z-\xi_l^{(k)}|\Br d\ml -\\ - \frac 1p\int_{\ql} \log \frac{|1-\bar z
\zeta|^p}{\prod_{k=1}^p |1-\bar z\xi_l^{(k)}|}d\ml =J_3+J_4.
\end{gather*}
As in the proof of the proposition (see the estimate of $I_2$), one can show that
 Since $|\frac1 {\bar z} -\zeta| \asymp 1- |z|\asymp |\frac1 {\bar z}
 -\xi_l^{(j)}|$. Hence,
we have $J_4=O(1)$.

Let  $\mu_z(t)=\mu(\ov{U(z,t)})$. Further, using
(\ref{e:zmz}),
\begin{gather} \nonumber
J_3=\int\limits_{\ql\setminus U(z,b(|z|))} \log|z-\zeta| d\ml+\int\limits_{U(z,
b(|z|))}  \log|z-\zeta| \, d\ml - \\ - p\log b(|z|) +O( 1)
=\mu^{(l)}(\ql\setminus U(z,b(|z|))\log b(|z|) +O(1) \nonumber +
\\
+\int_0^{b(|z|)} \log t \, d\mu_z^{(l)} (t) -p\log b(|z|)
=\mu^{(l)} (\ql \setminus U(z,b(|z|))\log b(|z|) + \nonumber\\ +O(1) +
\mu^{(l)} (U(z,b(|z|))\log b(|z|) - \int_0^{b(|z|)} \frac {\mu_z^{(l)}(t)}t
dt\nonumber
- \\ -p\log
b(|z|)
=-\int_0^{b(|z|)} \frac{\mu_z^{(l)} (t)}t dt+O(1)= O(1) \label{e:reg}
\end{gather}
by the regularity of $\mu_u$ w.r.t $b(t)$. Together with
(\ref{e:ej0}) it yields \bey \label{e:e0} \sum_{l\in \L^0}
|j_l(z)|=O(1), \quad z\not\in E_\ve. \eey

Now we estimate $\sum_{l\in \L^+} j_l(z)$.
Integration by parts gives us
\bey\label{e:ip}
L(\zeta)-L(\xi_l^{(1)})=\sum_{k=1}^m \frac 1{k!} L^{(k)}
(\xi_l^{(1)}) (\zeta -\xi_l^{(1)})^k + \frac 1{m!} \int_{\xi_l^{(1)}}^\zeta
L^{(m+1)} (s) (\zeta-s)^m \, ds ,
\eey
where $L(\zeta)=\log  \frac{z-\zeta}{1-\bar z\zeta} $,
\bey \label{e:ld}
|L^{(k)}(\zeta)| \le \frac {2(k-1)!}{|z-\zeta|^k}
\eey
The definition of $\xi_l^{(k)}$, $1\le k\le p$ allows us to cancel
the
first $p$  moments. Therefore, similarly to (\ref{e3.6}) and (\ref{e:dlfar})
 we have
\bey \label{e:jlfar}
|j_l(z)| \le C_{28} d_l^{p+1} \max_{s\in \ov{U(\ql, K_1(p)d_l)}}
|s-z|^{-p-1}.
\eey
Then
 $(z\in A_m)$
 \begin{gather*}
\sum_{l\in \L^+} |j_l(z)| \le C_{28}  \sum_{\begin{substack} {l\in \L^+
}\end{substack}}
\frac
{d_l^{p+1}}{|z- \xi_l^{(1)}|^{p+1}}
\le C_{29}  \sum_{l\in \L^+} d_l^{p-1}
\int_{\ql} \frac{dm(z)}{|z-\zeta|^{p+1}} \le \\ \le C_{30}(N,p,q)\sum_{n\le m-12}
b^{p-1}(R_n) \int_{\bar A_n} \frac{dm(z)}{|z-\zeta|^{p+1}} \le
 C_{31} \int\limits_{|\zeta|\le R_{m-12}} \frac {b^{p-1}(|\zeta|)
dm(\zeta)}{|z-\zeta|^{p+1}} \le  \\ \le C_{32} \int_0^{R_{m-12}} \frac
{b^{p-1}(\rho)}{(|z|-\rho)^p} d\rho \le C_{33} \int _0^1
\frac{b(\rho)^{p-1}}{(1-\rho)^p} \, d\rho<+\infty.
 \end{gather*}
Using the latter inequality, (\ref{e:e0}) and (\ref{e:ejm}) we
obtain  $|V(z)|=O(1)$ for $z\not \in E_\ve$.

The construction of $f$ is similar to that of Theorem \ref{t:1}.
It remains to prove (\ref{e:bo_ab}) for $z\in E_\ve$.

By (\ref{e:ej0}) it is sufficient to consider $l$ with $\ql \cap
U(z,3K_1(p)d_l)\ne \varnothing$.
For all sufficiently large $l\in \L^0$ we have
\begin{gather} \nonumber
\Bm \int_{\ql} \log|z-\zeta| d\ml\Bm \le
\int\limits_{U(z,4K_1(p)d_l)} \log \frac 1{|z-\zeta|} d\ml \le \\
\le \int\limits_0^{4K_1(p)d_l} \log \frac 1t d\mu_z^{(l)}(t) =\log
\frac 1{4K_1(p)d_l} \mu_z^{(l)}(4K_1(p)d_l)+
\int\limits_0^{4K_1(p)d_l} \frac{\mu_z^{(l)}(t)}t dt= O(1).
\label{e:3.20}
\end{gather}
Then we have \ben \log|f(z)| -u(z)=O(1) + \sum_{l\in \L^0} \Bl
\sum_{k=1}^p \log|z-\xi_l^{(k)}| -\int_{\ql} \log |z-\zeta|d\ml
\Br <C, \een because $|z-\xi_l^{(j)}|=O(b(|z|))<1$ for $l\ge l_0$
and (\ref{e:bo_ab}) is proved.

Finally, in order to prove (\ref{e:tt}) we note that for $z\in
E_\ve$ in view of (\ref{e:ejm}), (\ref{e:ej0}), (\ref{e:3.20})
 we have
$$\log|f(z)|- u(z)=\sum_{j=1}^m \log |z-\zeta_j|+O(1)$$ where
$\zeta_j\in Z_f$, and $m$ are uniformly bounded. Then $T(r,u-\log|f|)$
is bounded, and consequently $$ T(r,u)=T(r,\log|f|)+
T(r,u-\log|f|)+O(1)=T(r,\log|f|)+O(1).$$
\end{proof}

\begin{proof}[Proof of Theorem \ref{t:2}]
Let $\mu_j=\mu_u\Bigl|_{[\gamma_j]}$. By the assumptions of the
theorem we have $\mu_u=\sum_{j=1}^m \mu_j$. We can write
$u=\sum_{j=1}^m u_j$, where $u_j\in \SH(\mathbb{D})$, and
$\mu_{u_j} =\mu_j$. Therefore, it is sufficient to approximate
each $u_j$, $1\le j\le m$, separately.

We write $R(r)=(1-r)^{-1}$ and $W(R)=R^{\rho(R)}$. Then
$\mu_j(U(0,r))=\Delta_j W(R(r))$. Put $b(t)=(1-t)/W(R(r))$. Then
condition (\ref{e:bco}) is satisfied. We are going to prove that
$\mu_j$ admits a partition of slow variation and is locally
regular w.r.t. $b(t)$. We define a sequence $(r_n)$ from the
relation $\Delta_j W(R(r_n))=2n$, $n\in \N$. Then using the theorem on the inverse function,
 and properties
of the proximate order \cite[Ch.1, \S 12]{Le} we have ($r'\in (r_n, r_{n+1})$)
\begin{gather*}
r_{n+1}-r_n= \frac{2(1-r')^2 }{\Delta_j W'(R(r'))}=\frac
{(2+o(1))R(r')(1-r)'^2}{\Delta_j\sigma
W(R(r'))}=\frac{2+o(1)}{\Delta_j\sigma} b(r') \asymp b(r_n).
\end{gather*}
Let $Q^{(n)} =\{z: r_n\le |z|\le r_{n+1}, \vfi_n^- \le \theta \le
\vfi_{n}^+\}$ where $\vfi_n^-=\theta _j(r_n)-K(r_{n+1}-r_n)$,
$\vfi_{n}^+ =\theta _j(r_n)+K(r_{n+1}-r_n)$. Since
$|\theta'_j(t)|\le K$, we have $\theta_j(r)
\in[\vfi_{n}^-,\vfi_{n}^+]$, $r_{n}\le r\le r_{n+1}$. Let
$\mu^{(n)} =\mu_j\Bigl|_{Q^{(n)}}$. Then, by the definition of
$r_n$, $\mu^{(n)}(Q^{(n)})=2$. Therefore conditions i) and  iv) in
the definition of a partition of slow growth are satisfied.
Condition ii) is trivial. Since $\diam Q^{(n)}\asymp b(r_n)\asymp
(1-r_n)^{1+\si(r_n)} $, $\si>0$, conditions iii) and v) are valid.
Therefore, $\mu$ admits a partition of slow growth w.r.t. $b$,
$N=p=2$.

Finally, we check the local regularity of $\mu_j$ w.r.t. $b(t)$.
For $|z|=r$, $\rho\le b(r)$ we have
\begin{gather*} \mu_j(U(z,\rho)) \le
\Delta_j W(R(r+\rho))-\Delta_j W(R(r-\rho))=  W'(R(r^*))\frac
{2\rho\Delta_j}{(1-r^*)^2} =
\\
=
(2+o(1))\Delta_j\sigma\rho \frac{W(R(r^*))}{1-r^*} \le \frac
{3\sigma\rho\Delta_j}{b(r)}.
\end{gather*}
Then $\int_0^{b(r)} \frac{\mu(U(z,\rho))}\rho d\rho\le 3\si
\Delta_j$ as required.

Applying Theorem \ref{t:3} we obtain (\ref{e:uni1}),
(\ref{e:uni11}), and (\ref{e:tob}) for some analytic function
$f_j$ in $\mathbb{D}$.

Finally, we define $f=\prod_{j=1}^m f_j$.

The theorem is proved.
\end{proof}

I would like to thank  Professor O. Skaskiv who read the paper
and made valuable suggestion as well as other participants of
the Lviv  seminar on the theory of analytic
functions for valuable comments which contribute to the improvement the
initial version of the paper. I wish also to thank the anonymous referee for a lot
useful comments and corrections.

\renewcommand{\refname}{References}

\end{document}